\documentclass[a4paper,12pt]{article}
\usepackage[T1]{fontenc}
\usepackage[latin1]{inputenc}
\usepackage[english]{babel}
\usepackage{hyperref}
\usepackage[dvips]{graphics}
\usepackage{amsmath}
\usepackage{amssymb}
\usepackage{amsopn}
\usepackage{amsbsy}
\usepackage{amstext}
\usepackage{latexsym}
\usepackage{amsfonts}
\usepackage{array}
\usepackage{epsfig}
\usepackage{color}
\usepackage{mathrsfs} %%% \mathscr
\title{Algebraic and topological structures on the set of mean
functions and generalization of the AGM mean}
\author{\sc Bakir FARHI}
\date{}

\newtheorem{thm}{Theorem}[section]
\newtheorem{prop}[thm]{Proposition}

\newtheorem{rmq}[thm]{Remark}
\newtheorem{coll}[thm]{Corollary}

\newtheorem{defi}[thm]{Definition}

\let\epsilon=\varepsilon
\def\EMts{\mspace{.3mu}}

\def\G{{\rm G}}

\def\d{{\underline{\bf{d}}}}
\def\nb#1{{\left\vert{\EMts\EMts #1 \EMts\EMts}\right\vert}}

\def\EMdash{\leavevmode\hbox to 7.5mm{\vrule height .63ex depth -.59ex
    width 5.4mm\hfill}}
\def\D{\mathscr{D}}
\def\R{\mathbb{R}}
\def\A{{\rm A}}
\def\G{{\rm G}}
\def\H{{\rm H}}
\def\AGM{{\rm AGM}}
\def\M{{\mathcal{M}}}
\def\AS{{\mathcal{A}}}
\def\d{{\mathbf{\rm d}}}

\begin{document}
\maketitle \vspace{-7cm}
\begin{flushleft}
To appear
\end{flushleft}~\vspace{4.5cm}

\begin{abstract}
In this paper, we present new structures and results on the set
$\M_\D$ of mean functions on a given symmetric domain $\D$ of
$\mathbb{R}^2$. First, we construct on $\M_\D$ a structure of
abelian group in which the neutral element is simply the {\it
Arithmetic} mean; then we study some symmetries in that group.
Next, we construct on $\M_\D$ a structure of metric space under
which $\M_\D$ is nothing else the closed ball with center the {\it
Arithmetic} mean and radius $1/2$. We show in particular that the
{\it Geometric} and {\it Harmonic} means lie in the border of
$\M_\D$. Finally, we give two important theorems generalizing the
construction of the $\AGM$ mean. Roughly speaking, those theorems
show that for any two given means $M_1$ and $M_2$, which satisfy
some regular conditions, there exists a unique mean $M$ satisfying
the functional equation: $M(M_1 , M_2) = M$.
\end{abstract}
\vspace{1mm}
\section{Introduction}~

Let $\D$ be a nonempty symmetric domain of $\R^2$. A {\it mean}
function (or simply a {\it mean}) on $\D$ is a function $M : \D
\rightarrow \R$ satisfying the three following axioms:
\begin{description}
\item[i)] $M$ is symmetric, that is for all $(x , y) \in \D$:
$$M(x , y) ~=~ M(y , x) .$$
\item[ii)] For all $(x , y) \in \D$, we have:
$$\min(x , y) ~\leq~ M(x , y) ~\leq~ \max(x , y) .$$
\item[iii)] For all $(x , y) \in \D$, we have:
$$M(x , y) ~=~ x \Longrightarrow x ~=~ y .$$
\end{description}
Let us just remark that because of the axiom ii), the implication
of the axiom iii) is actually even an equivalence.\\
Among the most known examples of mean functions, we cite:\\
$\bullet$ The arithmetic mean, noted $\A$ and defined on $\R^2$
by:
$$\A(x , y) ~=~ \frac{x + y}{2} .$$
$\bullet$ The geometric mean, noted $\G$ and defined on ${(0 , +
\infty)}^2$ by:
$$\G(x , y) = \sqrt{x y} .$$
$\bullet$ The harmonic mean, noted $\H$ and defined on ${(0 , +
\infty)}^2$ by:
$$\H(x , y) = \frac{2 x y}{x + y} .$$
$\bullet$ The Gauss arithmetic-geometric mean, noted $\AGM$ and
defined on ${(0 , + \infty)}^2$ by the following process:\\
Given $x , y$ positive real numbers, $\AGM(x , y)$ is the common
limit of the two adjacent sequences ${(x_n)}_{n \in \mathbb{N}}$
and ${(y_n)}_{n \in \mathbb{N}}$ defined by:
$$\left\{\begin{array}{l} x_0 = x ~,~ y_0 = y \\ x_{n + 1} = \frac{x_n + y_n}{2} ~~~~~~(\forall n \in \mathbb{N})
\\ y_{n + 1} = \sqrt{x_n y_n} ~~~~~~ (\forall n \in \mathbb{N})\end{array}\right..$$

For a more profound survey on the mean functions, we refer to the
chapter 8 of the book \cite{bor} in which the mean $\AGM$ takes
the principal place. However, there are some differences between
that reference and the present paper. Indeed, in \cite{bor}, only
the axiom ii) is considered for defining a mean function; the
axiom iii) is added for obtaining the so called {\it strict} mean
while the axiom of symmetry ii) is not taken into account. In this
paper, we shall see that the three axioms i), ii) and iii) are
both necessary and sufficient to define a {\it good} mean or a
{\it good} set of mean functions on a given domain. In particular,
the axiom ii), which is excluded in the book \cite{bor}, is
necessary for the foundation of our algebraic and topological
structures announced in the title of this paper (see sections 2
and 3).~\vspace{1mm}

\noindent{\bf Remark about the axiom iii).} The axiom iii) permits
to avoid functions as $(x , y) \mapsto \min(x , y)$ and $(x , y)
\mapsto \max(x , y)$ which are not {\it true} means although they
satisfy the two axioms i) and ii). But beyond this simple
constatation, the axiom iii) will play a vital role throughout
this paper and especially in the generalization of the
construction of the arithmetic-geometric mean established by
Theorem \ref{t2}.

Given a nonempty symmetric domain $\D$ of $\R^2$, we denote by
$\M_{\D}$ the set of mean functions on $\D$. The purpose of this
paper is on the one hand to establish important algebraic and
topological structures on $\M_{\D}$ and to study some of their
properties and on the other hand to generalize in a natural way
the arithmetic-geometric mean $\AGM$. The article is composed of
three sections:\\
In the first section, we define on $\M_\D$ a structure of abelian
group in which the neutral element is simply the arithmetic mean
on $\D$. The study of this group reveals us that the arithmetic,
geometric and harmonic means lie in a very particular class of
mean functions that we call {\it the normal mean functions}. We
then study the symmetries on $\M_\D$ and we discover that the
symmetries with respect to one of the three means $\A$, $\G$ and
$\H$ oddly coincides with another type of symmetry (with respect
to the same means) we introduce and call {\it the functional
symmetry}. The problem of describing the set of the all means
realizing that curious coincidence is still open.\\
In the second section, we define on $\M_\D$ a structure of metric
space which turns out to be a closed ball with center $\A$ and
radius $1/2$. We then use the group structure (introduced in the
first section) to calculate the distance between two arbitrary
means on $\D$; this permits us in particular to establish a simple
characterization of the border of $\M_\D$.\\
In the third section, we introduce the concept of {\it functional
middle} of two mean functions on $\D$ which generalize in a
natural way the arithmetic-geometric mean, so that the latter is
the functional middle of the arithmetic mean and the geometric
mean. We establish two theorems, each provides a sufficiently
condition for the existence and the uniqueness of the functional
middle of two given means. The first one uses the topological
structure of $\M_\D$ by imposing on the two means in question the
condition that they are not diametrically opposed. The second one
(more important) imposes on the two means in question to be just
continuous on $\D$. In the proof of the latter one, the axiome
iii) plays an extremely vital role.

\section{An abelian group structure on $\M_\D$}~

For the following, given $\D$ a nonempty symmetric domain of
$\R^2$, we call $\AS_\D$ the set of asymmetric maps on $\D$; that
is maps $f : \D \rightarrow \R$, satisfying:
$$f(x , y) ~=~ - f(y , x) ~~~~~~ (\forall (x , y) \in \D) .$$
It is clair that $(\AS_\D , +)$ (where $+$ is the usual addition
law of the maps from $\D$ into $\R$) is an abelian group with
neutral element the null map from $\D$ into $\R$.\\
Now, consider $\varphi : \M_\D \rightarrow \R^{\D}$ the map
defined by:
$$\forall M \in \M_\D ~,~ \forall (x , y) \in \D :~~~~ \varphi(M)(x , y) :=
\begin{cases} \log\left(- \frac{M(x , y) - x}{M(x , y) - y}\right) & \text{if $x \neq y$} \\ 0 & \text{if x = y}\end{cases} .$$
The axioms i), ii) and iii) (verified by $M$, as a mean function)
insure the well-definition of $\varphi$, that is they insure that
the quantity $- \frac{M(x , y) - x}{M(x , y) - y}$ (for $x \neq
y$) is well-defined and positive.\\
The axiom of symmetry i) shows in addition that for all $M \in
\M_\D$, we have $\varphi(M) \in \AS_D$. Indeed, for all $M \in
\M_\D$ and for all $(x , y) \in \D$ such that $x \neq y$, we have:
\begin{eqnarray*}
\varphi(M)(x , y) & = & \log\left(- \frac{M(x , y) - x}{M(x , y) -
y}\right) \\
& = & \log\left(- \frac{M(y , x) - x}{M(y , x) - y}\right) ~~~~~~
\text{(according to the axiom i))} \\
& = & - \log\left(- \frac{M(y , x) - y}{M(y , x) - x}\right) \\
& = & - \varphi(M)(y , x) .
\end{eqnarray*}
Since we have in addition $\varphi(M)(x , x) = 0$ $(\forall (x ,
x) \in \D)$, then it follows that $\varphi(M)$ is effectively an
asymmetric map on $\D$, as claimed.\\
Conversely, if $f$ is an asymmetric map on $\D$, we easily verify
that $M : \D \rightarrow \R$, defined by:
$$M(x , y) := \frac{x + y e^{f(x , y)}}{e^{f(x , y)} + 1} ~~~~~~ (\forall (x , y) \in \D)$$
is a mean function on $\D$ and that we have $\varphi(M) = f$.\\
So, it follows that the map $\varphi$ constitutes a bijection from
$\M_\D$ to $\AS_\D$ and that its inverse map $\varphi^{-1}$ is
given by:
\begin{equation}\label{eq0}
\forall f \in \AS_D ~,~ \forall (x , y) \in \D : ~~~~
\varphi^{-1}(f)(x , y) = \frac{x + y e^{f(x , y)}}{e^{f(x , y)} +
1} .
\end{equation}
We thus can transport, by $\varphi$, the abelian group structure
$(\AS_\D , +)$ on $\M_\D$. This consists to define on $\M_\D$ the
following internal composition law $*$:
$$\forall M_1 , M_2 \in \M_\D : ~~~~ M_1 * M_2 = \varphi^{-1}\left(\varphi(M_1) + \varphi(M_2)\right) .$$
So $(\M_\D , *)$ is an abelian group and $\varphi$ is a group isomorphism from $(\M_\D , *)$ to $(\AS_D , +)$.\\
Furthermore, since the null map on $\D$ is the neutral element of
the group $(\AS_\D , +)$ and that $\varphi^{-1}(0) = \A$, then we
deduce that the arithmetic mean
$\A$ on $\D$ is the neutral element of the group $(\M_\D , *)$.\\
By calculating explicitly $M_1 * M_2$ (for $M_1 , M_2 \in \M_\D$),
we obtain the following theorem:
\begin{thm}
Let $\D$ be a nonempty symmetric domain of $\R^2$. Then, the law
$*$ on $\M_\D$ defined by: $\forall M_1 , M_2 \in \M_\D ~,~
\forall (x , y) \in \D$:
$$(M_1 * M_2)(x , y) := \begin{cases}
\frac{x (M_1(x , y) - x) (M_2(x , y) - y) + y (M_1(x , y) - x)
(M_2(x , y) - y)}{(M_1(x , y) - x) (M_2(x , y) - y) + (M_1(x , y)
- x) (M_2(x , y) - y)} & \text{if $x \neq y$} \\ x & \text{if x =
y}\end{cases}$$ is an internal composition law on $\M_\D$ and
$(\M_\D , *)$ is an abelian group with neutral element the
arithmetic mean $\A$ on $\D$.\\
In addition, the map $\varphi : \M_\D \rightarrow \AS_D$ defined
by:
$$\forall M \in \M_\D ~,~ \forall (x , y) \in \D :~~~~ \varphi(M)(x , y) :=
\begin{cases} \log\left(- \frac{M(x , y) - x}{M(x , y) - y}\right) & \text{if $x \neq y$} \\ 0 & \text{if x = y}\end{cases}$$
constitutes a group isomorphism from $(\M_\D , *)$ to $(\AS_D ,
+)$.~$\hfill\blacksquare$
\end{thm}
Now, let us calculate the images of the geometric and harmonic
means by the isomorphism $\varphi$ (with $\D = {(0 , +
\infty)}^2$). For all $(x , y) \in \D$, $x \neq y$, we have:
\begin{eqnarray*}
\varphi(\G)(x , y) & = & \log\left(- \frac{\G(x , y) - x}{\G(x ,
y) - y}\right) \\
& = & \log\left(- \frac{\sqrt{x y} - x}{\sqrt{x y} - y}\right) \\
& = & \log\left(- \frac{\sqrt{x} (\sqrt{y} - \sqrt{x})}{\sqrt{y}
(\sqrt{x} - \sqrt{y})}\right) \\
& = & \log\sqrt{x} - \log\sqrt{y} \\
& = & \frac{1}{2} \log{x} - \frac{1}{2} \log{y} .
\end{eqnarray*}
Since, in addition, $\varphi(\G)(x , x) = 0$ $(\forall x \in (0 ,
+\infty))$, it follows that:
\begin{equation}\label{eq1}
\varphi(\G)(x , y) ~=~ \frac{1}{2} \log{x} - \frac{1}{2} \log{y}
~~~~~~ (\forall (x , y) \in \D) .
\end{equation}
Similarly, for all $(x , y) \in \D$, $x \neq y$, we have:
\begin{eqnarray*}
\varphi(\H)(x , y) & = & \log\left(- \frac{\H(x , y) - x}{\H(x ,
y) - y}\right) \\
& = & \log\left(- \frac{\frac{2 x y}{x + y} - x}{\frac{2 x y}{x +
y} - y}\right) \\
& = & \log\left(\frac{x}{y}\right) \\
& = & \log{x} - \log{y} .
\end{eqnarray*}
Since we have in addition $\varphi(\H)(x , x) = 0$ $(\forall x \in
(0 , + \infty))$, it follows that:
\begin{equation}\label{eq2}
\varphi(\H)(x , y) ~=~ \log{x} - \log{y} ~~~~~~ (\forall (x , y)
\in \D) .
\end{equation}
From (\ref{eq1}) and (\ref{eq2}), we constat that $\varphi(\G)$
and $\varphi(\H)$ (with also $\varphi(\A)$) have a particular form
of asymmetric maps: each of them can be written
as $h(x) - h(y)$, where $h$ is a real function of one variable.\\
Conversely, given a nonempty subset $I$ of $\R$ and a map $h : I
\rightarrow \R$, it is clair that the map $f : I^2 \rightarrow
\R$, defined by:
$$f(x , y) ~=~ h(x) - h(y) ~~~~~~ (\forall x , y \in I)$$
is an asymmetric map on $I^2$. Consequently, $\varphi^{-1}(f)$
will give a mean on $I^2$. Let us explicit the expression of that
mean. For all $(x , y) \in I^2$, we have:
\begin{eqnarray*}
\varphi^{-1}(f)(x , y) & = & \frac{x + y e^{f(x , y)}}{e^{f(x , y)} + 1} ~~~~~~ \text{(according to (\ref{eq0}))}\\
& = & \frac{x + y e^{h(x) - h(y)}}{e^{h(x) - h(y)} + 1} \\
& = & \frac{x e^{-h(x)} + y e^{-h(y)}}{e^{-h(x)} + e^{-h(y)}} \\
& = & \frac{x P(x) + y P(y)}{P(x) + P(y)} ,
\end{eqnarray*}
with $P(t) = e^{-h(t)}$ $(\forall t \in I)$. We remark that the
only particularity of $P$ is that it is a positive function on
$I$. This leads us to include the most known three means
(arithmetic, geometric and harmonic means) in an important class
of mean functions which we define in what follows:
\begin{defi}{\rm
Let $I$ be a nonempty subset of $\R$. We call {\it normal} mean
function on $I^2$ any function $M : I^2 \rightarrow \R$, which can
be written as:
$$M(x , y) ~=~ \frac{x P(x) + y P(y)}{P(x) + P(y)} ~~~~~~ (\forall x , y \in I) ,$$
with $P : I \rightarrow \R$ is a positive function on $I$.\\
We call $P$ the {\it weight} function associated to $M$.
}\end{defi}

\noindent{\bf Remark.} The weight function associated to a normal
mean function is defined up to a multiplicative positive constant.
The weight functions associated to the three means $\A$, $\G$ and
$\H$ are respectively $P_\A(x) = 1$ $(\forall x \in \R)$, $P_\G(x)
= 1/\sqrt{x}$ $(\forall x > 0)$ and $P_\H(x) = 1/x$ $(\forall x >
0)$.

Before continuing our study on the group structure defined above
on $\M_\D$, we would like to stress a quit interesting property
concerning the comparison (in the sense of the usual order on
$\R$) between two normal mean functions.\\
To compare between two normal mean functions, there is a simple
and practical criterium which uses their associated weight
functions. We have the following
\begin{prop}\label{prop1}
Let $I$ be a nonempty interval of $\R$ and $M_1$ and $M_2$ be two
normal mean functions on $I^2$ with weight functions $P_1$ and
$P_2$ respectively. Then, the two following properties are
equivalent:
\begin{description}
\item[1)] $\forall (x , y) \in  I^2$: $M_1(x , y) \leq M_2(x ,
y)$.
 \item[2)] The function $\frac{P_1}{P_2}$ is non-increasing on $I$.
\end{description}
The same holds for the following two properties:
\begin{description}
\item[3)] $\forall (x , y) \in  I^2$, $x \neq y$: $M_1(x , y) <
M_2(x , y)$.
 \item[4)] The function $\frac{P_1}{P_2}$ is decreasing on $I$.
\end{description}
\end{prop}

\noindent{\bf Proof.} We only prove the equivalence of the two
properties 1) and 2). The prove of the equivalence of the two
properties 3) and 4) is similar. For all $(x , y) \in I^2$, we
have:
$$
M_1(x , y) \leq M_2(x , y) \Leftrightarrow \frac{x P_1(x) + y
P_1(y)}{P_1(x) + P_1(y)} \leq \frac{x P_2(x) + y P_2(y)}{P_2(x) +
P_2(y)} ~~~~~~~~~~~~~~~~~~~~~~~~~~~~~~~~~~~~~~~~$$ ~\vspace{-7mm}
\begin{eqnarray*}
& \Leftrightarrow & \!\!\!\!\left(x P_1(x) + y
P_1(y)\right)\left(P_2(x) + P_2(y)\right) \leq \left(x P_2(x) + y
P_2(y)\right)\left(P_1(x) + P_1(y)\right) \\
& \Leftrightarrow & \!\!\!\! (x - y)\left(P_1(x) P_2(y) - P_1(y)
P_2(x)\right) \leq 0 \\
& \Leftrightarrow & \!\!\!\! (x - y)\left(\frac{P_1}{P_2}(x) -
\frac{P_1}{P_2}(y)\right) \leq 0 .
\end{eqnarray*}
So, the property 1) of the proposition is equivalent to the
property:
$$\forall (x , y) \in I^2 : (x - y)\left(\frac{P_1}{P_2}(x) -
\frac{P_1}{P_2}(y)\right) \leq 0 ,$$ which amounts to say that the
function $P_1/P_2$ is non-increasing on
$I$.\penalty-20\null\hfill$\blacksquare$\par\medbreak

From Proposition \ref{prop1}, we derive the following immediate
corollary:
\begin{coll}
Let $I$ be a nonempty interval of $\R$ and $M$ be a normal mean
function on $I^2$, with weight function $P$. Then
\begin{description}
\item[1)] $M$ is sub-arithmetic (that is $M$ verifies: $\forall (x
, y) \in I^2$: $M(x , y) \leq \A(x , y)$) if and only if $P$ is
non-increasing on $I$.
 \item[2)] $M$ is strictly sub-arithmetic (that is $M$ verifies:
$\forall (x , y) \in I^2$, $x \neq y$: $M(x , y) < \A(x , y)$) if
and only if $P$ is decreasing on $I$.
 \item[3)] $M$ is super-arithmetic (that is $M$ verifies:
$\forall (x , y) \in I^2$: $M(x , y) \geq \A(x , y)$) if and only
if $P$ is non-decreasing on $I$.
 \item[4)] $M$ is strictly super-arithmetic (that is $M$ verifies:
$\forall (x , y) \in I^2$, $x \neq y$: $M(x , y) > \A(x , y)$) if
and only if $P$ is increasing on $I$.
\penalty-20\null\hfill$\blacksquare$\par\medbreak
\end{description}
\end{coll}

\noindent{\bf An application of Proposition \ref{prop1}.} By
remembering that the respective weight functions associated to the
three normal means $\A$, $\G$ and $\H$ are $1$, $1/\sqrt{x}$ and
$1/x$, Proposition \ref{prop1} immediately shows that for all $x ,
y
> 0$, $x \neq y$, we have:
$$A(x , y) > G(x , y) > H(x , y) .$$

\subsection*{Study of some symmetries on the group $(\M_\D , *)$}~

Given a nonempty symmetric domain $\D$ of $\R^2$, we are
interested in what follows in the symmetric mean of a given mean
$M_1$ with respect to another given mean $M_0$ via the group
structure $(\M_\D , *)$. Denoting by $S_{M_0}$ the symmetry with
respect to a fixed mean $M_0$ in the group $(\M_\D , *)$, we have
by definition:
$$\forall M_1 , M_2 \in \M_\D : S_{M_0}(M_1) = M_2 \Longleftrightarrow M_1 * M_2 = M_0 * M_0 .$$
The explicit expression of $S_{M_0}(M_1)$ $(M_0 , M_1 \in \M_\D)$
is given by the following:
\begin{prop}\label{prop2}
Let $\D$ be a nonempty symmetric domain of $\R^2$ and $M_0$ and
$M_1$ be two mean functions on $\D$. Then we have:
$$S_{M_0}(M_1) = \frac{x (M_1 - x) (M_0 - y)^2 - y (M_0 - x)^2 (M_1 - y)}{(M_1 - x)(M_0 - y)^2 - (M_0 - x)^2 (M_1 - y)} ,$$
where, for simplicity, we have noted $M_0$ for $M_0(x , y)$, $M_1$
for $M_1(x , y)$ and $S_{M_0}(M_1)$ for $S_{M_0}(M_1)(x , y)$.
\end{prop}

\noindent{\bf Proof.} Of course, we use the group isomorphism
$\varphi$ introduced at the beginning of Section 2. Let $f_0 =
\varphi(M_0)$, $f_1 = \varphi(M_1)$, $M_2 = S_{M_0}(M_1)$ and $f_2
= \varphi(M_2)$. The equality $M_2 = S_{M_0}(M_1)$ amounts to $M_1
* M_2 = M_0 * M_0$. By applying $\varphi$ to the two sides of
the latter, we get $f_1 + f_2 = 2 f_0$, which gives $f_2 = 2 f_0 -
f_1$. It follows that for all $(x , y) \in \D$, we have:
\begin{eqnarray}
M_2(x , y) = \varphi^{-1}(f_2)(x , y) = \frac{x + y e^{f_2(x ,
y)}}{e^{f_2(x , y)} + 1} & = & \frac{x + y e^{2 f_0(x , y) - f_1(x
, y)}}{e^{2 f_0(x , y) - f_1(x , y)} + 1} \notag \\
& = & \frac{x e^{f_1(x , y)} + y e^{2 f_0(x , y)}}{e^{2 f_0(x ,
y)} + e^{f_1(x , y)}} \label{eq3}
\end{eqnarray}
Further, since $f_0 = \varphi(M_0)$ and $f_1 = \varphi(M_1)$, we
have:
$$e^{f_0(x , y)} = - \frac{M_0(x , y) - x}{M_0(x , y) - y} ~~~~\text{et}~~~~ e^{f_1(x , y)} = - \frac{M_1(x , y) - x}{M_1(x , y) - y} .$$
By inserting that two last equalities into (\ref{eq3}) and after
simplification and rearrangement, the identity of the proposition
follows.\penalty-20\null\hfill$\blacksquare$\par\medbreak

The following corollary gives us the expression of the symmetric
mean of a given mean with respect to one of the three means $\A$,
$\G$ and $\H$. Remark that the obtained expressions for the
symmetric means with respect to the arithmetic and geometric means
supports the intuition that we can have about them. This curious
fact shows the interest of the group structure defined on $\M_\D$.

\begin{coll}
Let $\D$ be a nonempty symmetric domain of $\R^2$ and $M$ a mean
function on $\D$. Then we have:
\begin{description}
\item[1)] $S_\A(M) = x + y - M = 2 \A - M$.
 \item[2)] $S_\G(M) = \frac{x y}{M} = \frac{\G^2}{M}$ ~~~~(by supposing $\D \subset {(0 , +
 \infty)}^2$).
 \item[3)] $S_\H(M) = \frac{x y M}{(x + y) M - x y} = \frac{\H M}{2 M -
 H}$ ~~~~(by supposing $\D \subset {(0 , + \infty)}^2$).
 \item[4)] $S_\H = S_\G \circ S_\A \circ S_\G$.
\end{description}
\end{coll}

\noindent{\bf Proof.} To obtain the formulas of the items 1), 2)
and 3), it suffices to apply the formula of Proposition
\ref{prop2} respectively for $M_0 = \A$, $M_0 = \G$ and $M_0 =
\H$. The formula of the item 4) is derived from those of the
previous one.\penalty-20\null\hfill$\blacksquare$\par\medbreak

Now, we are going to define another symmetry on the set $\M_\D$
($\D$ having a certain form) which is completely independent of
the group structure $(\M_\D , *)$. This new symmetry is defined by
solving a functional equation but it curiously coincides, in many
cases, with the symmetry defined above which is rather related to
the law of the group $(\M_\D , *)$.

\begin{defi}{\rm Let $I$ be a nonempty interval of $\R$, $\D = I^2$
and $M_0$, $M_1$ and $M_2$ be three mean functions on $\D$ such
that $M_1$ and $M_2$ take their values in $I$. We say that $M_2$
is the {\it functional symmetric mean} of $M_1$ with respect to
$M_0$ if the following functional equation is satisfied:
$$M_0(M_1(x , y) , M_2(x , y)) ~=~ M_0(x , y) ~~~~~~~~ (\forall (x , y) \in \D) .$$
Equivalently, we also say that $M_0$ is the {\it functional
middle} of $M_1$ and $M_2$. }\end{defi}

According to the axiom iii) verified by the mean functions, it is
immediate that if the functional symmetric of a mean with respect
to another mean exists then it is unique. This justifies the
following notation:~\vspace{1mm}

\noindent{\bf Notation.} Given $M_0$ and $M_1$ two mean functions
on $\D = I^2$ with values in $I$ (where $I$ is an interval of
$\R$), we denote by $\sigma_{M_0}(M_1)$ the functional symmetric
(if it exists) of $M_1$ with respect to $M_0$.~\vspace{1mm}

The following proposition gives us a sufficient condition for the
existence of the functional symmetric mean.
\begin{prop}\label{prop3}
Let $I$ be a nonempty interval of $\R$, $\D = I^2$ and $M_0$ be a
mean function on $\D$. For all $x \in I$, set:
$$I_x := \left\{M_0(x , t) ~|~ t \in I\right\}$$
and suppose that we have:
\begin{equation}\label{eq4}
\forall x , y \in I :~~~~ x \leq y \Longrightarrow I_x \subset I_y
.
\end{equation}
Suppose also that $M_0$ is increasing, that is $M_0$ is increasing
with respect to one (so to each) of its two variables.\\
Then, any mean function on $\D$ has a functional symmetric with
respect to $M_0$.
\end{prop}

\noindent{\bf Proof.} Let us fix a mean function $M_1$ on $\D$.
For all $(x , y) \in \D$, we are going to show the existence of a
unique real number, which we denote by $M_2(x , y)$, such that
\begin{equation}\label{eq5}
M_0(M_1(x , y) , M_2(x , y)) = M_0(x , y) .
\end{equation}
then, we will show that the obtained function $M_2$ (on $\D$) is
actually a mean on $\D$.\\
Let $x , y \in I$, fixed. consider the map:
$$
\begin{array}{rcl}
i : ~~ I & \longrightarrow & \R \\
t & \longmapsto & M_0(M_1(x , y) , t) .
\end{array}
$$
According to the hypothesis that $M_0$ is increasing, the map $i$
increases on $I$. It follows that $i$ is a bijection from $I$ into $i(I) = I_{M_1(x , y)}$.\\
Now, since $M_0(x , y) \in I_{\min(x , y)}$ (because $M_0(x , y)
\in I_x$ and $M_0(x , y) \in I_y$) and $I_{\min(x , y)} \subset
I_{M_1(x , y)}$ (according to the hypothesis (\ref{eq4})) then
$M_0(x , y) \in I_{M_1(x , y)}$. The bijectivity of $i$ from $I$
into $I_{M_1(x , y)}$ thus implies that there exists a unique $t_*
\in I$ satisfying $i(t_*) = M_0(x , y)$, so satisfying
$$M_0(M_1(x , y) , t_*) = M_0(x , y) .$$
It suffices to set $M_2(x , y) = t_*$ to establish (\ref{eq5})
for the couple $(x , y)$.\\
So the existence and the uniqueness of $M_2$ as a function on $\D$
satisfying the functional equation (\ref{eq5}) are confirmed and
it just remains to show that $M_2$ is a mean on $\D$.~\vspace{1mm}

\noindent $\bullet$ Let us show that $M_2$ is symmetric on $\D$:\\
Given $(x , y) \in \D$, $M_2(y , x)$ is (by definition) the unique
real number of $I$ satisfying
$$M_0(M_1(y , x) , M_2(y , x)) = M_0(y , x) .$$
But because $M_0$ and $M_1$ are symmetric on $\D$ (as mean
functions on $\D$), the last relation amounts to
$$M_0(M_1(x , y) , M_2(y , x)) = M_0(x , y) ,$$
which implies (according to the definition of $M_2(x , y)$) that:
$$M_2(x , y) = M_2(y , x) .$$
So the function $M_2$ is symmetric on $\D$ as we claimed it to be.

\noindent $\bullet$ Let us show that $M_2$ satisfies:
\begin{equation}\label{eq6}
\forall (x , y) \in \D : ~~~~ \min(x , y) \leq M_2(x , y) \leq
\max(x , y) .
\end{equation}
We argue by contradiction. Assume that there exists a couple $(x ,
y) \in \D$ for which (\ref{eq6}) is not valid. So, we
have:~\vspace{1mm}

\noindent\underline{Either $M_2(x , y) < \min(x , y)$:}\vspace{1mm}\\
In this case, since $M_1(x , y) \leq \max(x , y)$ (because $M_1$
is a mean on $\D$), then according to the hypothesis that $M_0$ is
increasing, we have:
$$M_0(M_1(x , y) , M_2(x , y)) < M_0(\max(x , y) , \min(x , y)) = M_0(x , y) ,$$
which contradicts the relation (\ref{eq5}) satisfied by $M_2$.
This first case is thus impossible.~\vspace{1mm}

\noindent\underline{Or $M_2(x , y) > \max(x ,
y)$:}\vspace{1mm}\\
In this case, since $M_1(x , y) \geq \min(x , y)$ (because $M_1$
is a mean on $\D$), then according to the hypothesis that $M_0$ is
increasing, we have:
$$M_0(M_1(x , y) , M_2(x , y)) > M_0(\min(x , y) , \max(x , y)) = M_0(x , y) ,$$
which again contradicts the relation (\ref{eq5}) satisfied by
$M_2$. This second case is thus also impossible. \\
The function $M_2$ thus satisfies the property (\ref{eq6}).

\noindent $\bullet$ Let us finally show that $M_2$ satisfies the
third axiom of mean functions, that is:
$$\forall (x , y) \in \D : ~~~~ M_2(x , y) = x \Longrightarrow x = y .$$
For all $(x , y) \in \D$, we have:
\begin{eqnarray*}
M_2(x , y) = x & \Longrightarrow & M_0(M_1(x , y) , M_2(x ,
y)) = M_0(M_1(x , y) , x) \\
& \Longrightarrow & M_0(x , y) = M_0(x , M_1(x , y)) \\
& \Longrightarrow & M_1(x , y) = y \\
& \Longrightarrow & x = y ,
\end{eqnarray*}
where in the second implication, we have used the relation
(\ref{eq5}) together with the symmetry of $M_0$; in the third
implication, we have used the increasing of $M_0$ and in the forth
implication, we have used the third axiom of mean functions for
the mean $M_1$.\\
The third axiom of mean functions is thus satisfied by $M_2$. \\
In conclusion, $M_2$ is a mean function on $\D$. This completes
the proof of the
proposition.\penalty-20\null\hfill$\blacksquare$\par\medbreak

We can see easily that the three means Arithmetic (on $\D =
\R^2$), Geometric and Harmonic (on $\D = {(0 , + \infty)}^2$)
satisfy the hypothesis of Proposition \ref{prop3}. So, we can
establish for any given mean (on a suitable domain $\D$) its
symmetric mean (in the functional sense) with respect to one of
the three means $\A$, $\G$ and $\H$. Precisely, We have the
following
\begin{prop}
Let $M$ be a mean function on a suitable symmetric domain $\D$ of
$\R^2$. Then the functional symmetric means of $M$ with respect to
the three means Arithmetic, Geometric and Harmonic are
respectively given by:
\begin{eqnarray*}
\sigma_A(M) & = & x + y - M \\
\sigma_G(M) & = & \frac{x y}{M} ~~~~~~~~~~~~~~~~~~~~~~~~ (\text{for $\D \subset {(0 , + \infty)}^2$}) \\
\sigma_H(M) & = & \frac{x y M}{(x + y) M - x y} ~~~~~~~~
(\text{for $\D \subset {(0 , + \infty)}^2$}) .
\end{eqnarray*}
So, those functional symmetric means coincide with the symmetric
means in the sense of the group law defined above on $\M_\D$.
\end{prop}

\noindent{\bf Proof.} Given $M$ a mean function on $\R^2$, its
functional symmetric mean with respect to the arithmetic mean $\A$
is defined by the functional equation:
$$\A(M , \sigma_\A(M)) = \A ,$$
which amounts to:
$$\frac{M + \sigma_\A(M)}{2} = \frac{x + y}{2} .$$
Hence:
$$\sigma_\A(M) = x + y - M ,$$
as claimed.\\
Similarly, given $M$ a mean function on ${(0 , + \infty)}^2$, its
functional symmetric means with respect to the two means $\G$ and
$\H$ are respectively defined by:
$$\G(M , \sigma_\G(M)) = \G ~~~~\text{and}~~~~ \H(M , \sigma_\H(M)) = \H ,$$
that is:
$$\sqrt{M \sigma_\G(M)} = \sqrt{x y} ~~~~\text{and}~~~~ \frac{2 M \sigma_H(M)}{M + \sigma_\H(M)} = \frac{2 x y}{x + y} .$$
This gives:
$$\sigma_\G(M) = \frac{x y}{M} ~~~~\text{and}~~~~ \sigma_\H(M) = \frac{x y M}{(x + y) M - x y} ,$$
as claimed. The proposition is
proved.\penalty-20\null\hfill$\blacksquare$\par\medbreak

The remarkable phenomenon of the coincidence of the two symmetries
defined on $\M_\D$ in the three particular cases corresponding to
the most known means $\A$, $\G$ and $\H$ leads us to formulate the
following important open question:

\noindent {\bf Open question.} What are all the mean functions $M$
on $D = {(0 , + \infty)}^2$ for which the two symmetries with
respect to $M$ (in the sense of the group law
introduced on $\M_\D$ and in the functional sense) coincide?\\

Mathematically speaking, we ask about the description of the set
$$\left\{M \in \M_{{(0 , + \infty)}^2} : S_M = \sigma_M\right\}$$
which contains (as proved above) at least the three means $\A$,
$\G$ and $\H$.~\vspace{1mm}

We end this section by giving an important example of functional
symmetry. We have the following

\begin{prop}
The two means $\A$ and $\G$ are symmetric in the functional sense
with respect to the $\AGM$ mean.
\end{prop}

\noindent{\bf Proof.} Given $x , y > 0$, each of the two real
numbers $\AGM(x , y)$ and $\AGM(\frac{x + y}{2} , \sqrt{x y})$ is
defined as the common limit of two adjacent sequences. But, it is
easy to see that the two adjacent sequences defining $\AGM(\frac{x
+ y}{2} , \sqrt{x y})$ are the shifted by one term of the two
adjacent sequences defining $\AGM(x , y)$. Consequently, we have:
$$\AGM\!\left(\frac{x + y}{2} , \sqrt{x y}\right) = \AGM(x , y) ,$$
that is:
$$\AGM(\A(x , y) , \G(x , y)) = \AGM(x , y) .$$
This shows that the means $\A$ and $\G$ are symmetric, in the
functional sense, with respect to the $\AGM$ mean. The proposition
is proved.\penalty-20\null\hfill$\blacksquare$\par\medbreak

\section{A metric space structure on $\M_\D$}~

Throughout this section, we fixe a nonempty symmetric domain $\D$
of $\R^2$. In the case where all the points of $\D$ have the form
$(x , x)$ $(x \in \R)$, the set $\M_\D$ is reduced to a unique
element and consequently there is a unique topology on $\M_\D$
which is trivial. So suppose that $\D$ contains at least a point
$(x_0 , y_0)$ of $\R^2$ such that $x_0 \neq y_0$. For all couple
$(M_1 , M_2)$ of mean functions on $\D$, define:
$$\d(M_1 , M_2) := \sup_{(x , y) \in \D , x \neq y} \nb{\frac{M_1(x , y) - M_2(x , y)}{x - y}} .$$
We have the following:
\begin{prop}\label{prop5}
The map $\d$ of $\M_\D^2$ into $[0 , + \infty]$ is a distance on
$\M_\D$. In addition, the metric space $(\M_\D , \d)$ is identic
to the closed ball with center $\A$ (the arithmetic mean) and
radius $\frac{1}{2}$.
\end{prop}

\noindent{\bf Proof.} First let us show that for all couple $(M_1
, M_2)$ of $\M_\D^2$, the nonnegative quantity $\d(M_1 , M_2)$ is
finite. Given $M_1 , M_2 \in \M_\D$, for all $(x , y) \in \D$, $x
\neq y$, the two real numbers $M_1(x , y)$ and $M_2(x , y)$ lie in
the same interval $[\min(x , y) , \max(x , y)]$, so we have:
$$\nb{M_1(x , y) - M_2(x , y)} \leq \max(x , y) - \min(x , y) = \nb{x - y} .$$
Hence
$$\sup_{(x , y) \in \D , x \neq y} \nb{\frac{M_1(x , y) - M_2(x , y)}{x - y}} \leq 1 ,$$
that is $\d(M_1 , M_2) \leq 1$.\\
This shows that $\d$ is actually a map from $\M_\D^2$ into $[0 ,
1]$. Further, since the three properties
\begin{quote}
$\bullet$ $\forall M_1 , M_2 \in \M_\D$: $\d(M_1 , M_2) = \d(M_2 ,
M_1)$ \\
$\bullet$ $\forall M_1 , M_2 \in \D$: $\d(M_1 , M_2) = 0
\Leftrightarrow M_1 = M_2$ \\
$\bullet$ $\forall M_1 , M_2 , M_3 \in \M_\D$: $\d(M_1 , M_3) \leq
\d(M_1 , M_2) + \d(M_2 , M_3)$
\end{quote}
are trivially satisfied by $\d$, then $\d$ is a distance on
$\M_\D$.

Now, given $M \in \M_\D$, let us show that $\d(M , \A) \leq
\frac{1}{2}$. For all couple $(x , y) \in \D$, $x \neq y$, the
real number $M(x , y)$ lies in the closed interval limited by $x$
and $y$, so we have:
\begin{eqnarray*}
\nb{M(x , y) - \A(x , y)} & \leq & \max\left(x - \A(x , y) , y
- \A(x , y)\right) \\
& = & \max\left(x - \frac{x + y}{2} , y - \frac{x + y}{2}\right)
\\
& = & \max(\frac{x - y}{2} , \frac{y - x}{2}) \\
& = & \frac{1}{2} \nb{x - y} .
\end{eqnarray*}
It follows that:
$$\sup_{(x , y) \in \D , x \neq y} \nb{\frac{M(x , y) - \A(x , y)}{x - y}} \leq \frac{1}{2} ,$$
that is $\d(M , \A) \leq \frac{1}{2}$, as required.\\
The metric space $(\M_\D , \d)$ is thus identic to the closed ball
with center $\A$ and radius $\frac{1}{2}$. This completes the
proof of the
proposition.\penalty-20\null\hfill$\blacksquare$\par\medbreak
\begin{rmq}\label{rmq1}
Given $M_1 , M_2 \in \M_\D$, since the map $(x , y) \mapsto
\frac{M_1(x , y) - M_2(x , y)}{x - y}$ is obviously asymmetric (on
the set $\{(x , y) \in \D : x \neq y\}$), we have also
$$\d(M_1 , M_2) = \sup_{(x , y) \in \D , x \neq y} \frac{M_1(x , y) - M_2(x , y)}{x - y} .$$
\end{rmq}

By using the group isomorphism $\varphi$ from $\M_\D$ into
$\AS_D$, defined at Section 2, we will establish in what follows a
practice formula to calculate the distance between two mean
functions on $\D$.
\begin{prop}\label{prop4}
Let $M_1$ and $M_2$ two mean functions on $\D$. Set $f_1 =
\varphi(M_1)$ and $f_2 = \varphi(M_2)$. Then we have:
$$\d(M_1 , M_2) = \sup_{(x , y) \in \D} \frac{e^{f_1} - e^{f_2}}{(e^{f_1} + 1) (e^{f_2} + 1)} = \sup_{(x , y) \in \D}
\left(\frac{1}{e^{f_1} + 1} - \frac{1}{e^{f_2} + 1}\right) .$$
\end{prop}

\noindent{\bf Proof.} According to the relation (\ref{eq0}) of
Section 2, we have for all $(x , y) \in \D$: $M_1(x , y) =
\varphi^{-1}(f_1)(x , y) = \frac{x + y e^{f_1(x , y)}}{e^{f_1(x ,
y)} + 1}$ and $M_2(x , y) = \varphi^{-1}(f_2)(x , y) = \frac{x + y
e^{f_2(x , y)}}{e^{f_2(x , y)} + 1}$. It follows that for all $(x
, y) \in \D$, we have:
\begin{eqnarray*}
M_1(x , y) - M_2(x , y) & = & \frac{x + y e^{f_1(x , y)}}{e^{f_1(x
, y)} + 1} - \frac{x + y e^{f_2(x , y)}}{e^{f_2(x , y)} + 1} \\
& = & (x - y) \frac{e^{f_2(x , y)} - e^{f_1(x ,
y)}}{\left(e^{f_1(x , y)} + 1\right)\left(e^{f_2(x , y)} +
1\right)} .
\end{eqnarray*}
Then, according to Remark \ref{rmq1} and the fact that the
functions $f_1$ and $f_2$ are zero at the points $(x , x)$ of $\D$
(because they are asymmetric on $\D$), we have:
$$\d(M_1 , M_2) = \sup_{(x , y) \in \D , x \neq y} \frac{M_1(x , y) - M_2(x , y)}{x - y} = \sup_\D \frac{e^{f_1}
- e^{f_2}}{(e^{f_1} + 1) (e^{f_2} + 1)} .$$ The proposition is
proved.\penalty-20\null\hfill$\blacksquare$\par\medbreak

Now, using Proposition \ref{prop4}, we will establish in what
follows a practice criterium which permits to locate easily any
mean function on $\D$ in the metric space $\M_\D$, seen as the
closed ball with center $\A$ and radius $\frac{1}{2}$.
\begin{coll}\label{coll1}
Let $M$ be a mean function on $\D$ and $f := \varphi(M)$, where
$\varphi$ is the group isomorphism defined at Section 2. Then,
setting $s := \sup_\D f \in [0 , + \infty]$, we have:
$$\d(M , \A) = \frac{1}{2} \cdot \frac{e^s - 1}{e^s + 1} $$
(We naturally suppose that $\frac{e^s - 1}{e^s + 1} = 1$
when $s = + \infty$). \\
In particular, the mean $M$ lies in the border of $\M_\D$ (that is
on the circle with center $\A$ and radius $\frac{1}{2}$) if and
only if $\sup_\D f = + \infty$.
\end{coll}

\noindent {\bf Proof.} Since the asymmetric function associated to
the arithmetic mean by the isomorphism $\varphi$ is the zero
function (i.e., $\varphi(\A) \equiv 0$), then according to
Proposition \ref{prop4}, we have:
$$\d(M , \A) = \sup_\D \frac{e^f - 1}{2 (e^f + 1)} = \frac{1}{2} \sup_\D \frac{e^f - 1}{e^f + 1} .$$
Next, since the function $x \mapsto \frac{e^x - 1}{e^x + 1}$ is
increasing on $\R$, then we have $\sup_\D \frac{e^f - 1}{e^f + 1}
= \frac{e^s - 1}{e^s + 1}$. The corollary
follows.\penalty-20\null\hfill$\blacksquare$\par\medbreak

\noindent{\bf Applications.}~
\begin{description}
\item[1)] The two means Geometric and Harmonic (on $\D = {(0 , +
\infty)}^2$) lie on the border of the metric space $\M_\D$.
Indeed, the asymmetric functions associated to $\G$ and $\H$ by
the isomorphism $\varphi$ are respectively: $\varphi(\G)(x , y) =
\frac{1}{2} (\log x - \log y)$ and $\varphi(\H)(x , y) = \log x -
\log y$ (see §2). Since we have clearly $\sup_\D \varphi(\G) =
\sup_\D \varphi(\H) = + \infty$, the second statement of Corollary
\ref{coll1} insures that $\G$ and $\H$ lie on the border of
$\M_\D$ as we claimed it to be.
 \item[2)] The distance between the two means Geometric and
 Harmonic (on $\D = {(0 , + \infty)}^2$) is somewhat more
 difficult to calculate. Using Proposition \ref{prop4}, we can show that:
 $$\d(\G , \H) = \sup_{t > 0} \frac{t^2 - t}{(t + 1)(t^2 + 1)} \simeq 0.15 .$$
Actually, we can show that this distance is an algebraic number
with degree 4; it is a root of the equation $x^4 + 10 x^3 + 3 x^2
- 14 x + 2 = 0$.
\end{description}

\section{Construction of a functional middle of two mean functions (generalization of the $\AGM$
mean)}~

Let $I$ be a nonempty subset of $\R$ and let $\D = \R^2$. The aim
of this section is to prove, under some {\it regular} conditions,
the existence and the uniqueness of the functional middle of two
given means $M_1$ and $M_2$ on $\D$; that is the existence and the
uniqueness of a new mean $M$ on $\D$ satisfying the functional
equation:
$$M(M_1 , M_2) = M .$$
In this context, we obtain two results which only differ in the
imposed condition on the two means $M_1$ and $M_2$. The first one
imposes to $M_1$ and $M_2$ the condition $\d(M_1 , M_2) \neq 1$
(where $\d$ is the distance on $\M_\D$ defined at Section 3) while
the second one simply imposes to $M_1$ and $M_2$ to be continuous
on $\D$ (by taking $I$ an interval of $\R$). Notice further that
our way of establishing the existence of the functional middle is
constructive and generalizes the idea of the $\AGM$ mean. Our
first result is the following:

\begin{thm}\label{t1}
Let $M_1$ and $M_2$ be two mean functions on $\D = I^2$, with
values in $I$ and such that $\d(M_1 , M_2) < 1$. Then there exists
a unique mean function $M$ on $\D$ satisfying the functional
equation:
$$M(M_1 , M_2) = M .$$
Besides, for all $(x , y) \in \D$, $M(x , y)$ is the common limit
of the two real sequences ${(x_n)}_n$ and ${(y_n)}_n$ defined as
follows:
$$\left\{\!\begin{array}{l}
x_0 = x ~,~ y_0 = y \\
x_{n + 1} = M_1(x_n , y_n)~~~~(\forall n \in \mathbb{N})\\
y_{n + 1} = M_2(x_n , y_n)~~~~(\forall n \in \mathbb{N})
\end{array}\right..$$
\end{thm}

\noindent{\bf Proof.} Let $k := \d(M_1 , M_2)$. By hypothesis, we
have $k < 1$. Let ${(x_n)}_n$ and ${(y_n)}_n$ be the two real
sequences introduced at the second statement of the theorem and
let ${(u_n)}_n$ and ${(v_n)}_n$ be the two real sequences (with
values in $I$) defined by:
$$u_n := \min(x_n , y_n) ~~\text{and}~~ v_n := \max(x_n , y_n) ~~~~~~ (\forall n \in \mathbb{N}) .$$
Let us show that ${(u_n)}_n$ and ${(v_n)}_n$ are adjacent
sequences. For all $n \in \mathbb{N}$, we have:
$$u_{n + 1} = \min(x_{n + 1} , y_{n + 1}) = \min(M_1(x_n , y_n) , M_2(x_n , y_n)) \geq \min(x_n , y_n) = u_n$$
(because $M_1(x_n , y_n) \geq \min(x_n , y_n)$ and $M_2(x_n , y_n)
\geq \min(x_n , y_n)$, since $M_1$ and $M_2$ are mean functions). This shows that ${(u_n)}_n$ is a non-decreasing sequence.\\
Similarly, for all $n \in \mathbb{N}$, we have:
$$v_{n + 1} = \max(x_{n + 1} , y_{n + 1}) = \max(M_1(x_n , y_n) , M_2(x_n , y_n)) \leq \max(x_n , y_n) = v_n$$
(because $M_1(x_n , y_n) \leq \max(x_n , y_n)$ and $M_2(x_n , y_n)
\leq \max(x_n , y_n)$, since $M_1$ et $M_2$ are mean functions). The sequence ${(v_n)}_n$ is then non-increasing.\\
Next, for all $n \in \mathbb{N}$, we have:
\begin{eqnarray*}
\nb{v_{n + 1} - u_{n + 1}} & = & \nb{\max(x_{n + 1} , y_{n + 1}) -
\min(x_{n + 1} , y_{n + 1})} \\
& = & \nb{x_{n + 1} - y_{n + 1}} \\
& = & \nb{M_1(x_n , y_n) - M_2(x_n , y_n)} \\
& \leq & k \nb{x_n - y_n} ~~~~~~~~\text{(by definition of $k$)} \\
& = & k \nb{v_n - u_n} .
\end{eqnarray*}
By induction on $n$, we get:
$$\nb{v_n - u_n} \leq k^n \nb{v_0 - u_0} ~~~~ (\forall n \in \mathbb{N}) .$$
It follows (since $k \in [0 , 1)$) that $(v_n - u_n)$ tends to
$0$ as $n$ tends to infinity.\\
The two sequences ${(u_n)}_n$ and ${(v_n)}_n$ are thus adjacent,
as claimed. Consequently, ${(u_n)}_n$ and ${(v_n)}_n$ are
convergent and have the same limit. In addition, since we have
clearly:
$$u_n \leq x_n \leq v_n ~~\text{and}~~ u_n \leq y_n \leq v_n ~~~~ (\forall n \in \mathbb{N}) ,$$
then ${(x_n)}_n$ and ${(y_n)}_n$ are also convergent and have the
same limit which coincides with the common limit of the two
sequences ${(u_n)}_n$ and ${(v_n)}_n$. In the sequel, let $M(x ,
y)$ denote this limit.

Now let us show that the map $M : \D \rightarrow \R$, just
defined, is a mean function on $\D$ and satisfies $M(M_1 , M_2) =
M$. First let us show that $M$ satisfies the three axioms of a
mean function.~\vspace{1mm}\\
{\bf i)} Given $(x , y) \in \D$, by changing $(x , y)$ by $(y ,
x)$ in the definition of the sequences ${(x_n)}_n$ and
${(y_n)}_n$, those lasts still unchanged except their first terms
(since $M_1$ and $M_2$ are symmetric). So, we have:
$$M(x , y) = M(y , x) ~~~~ (\forall (x , y) \in \D) .$$
In other words, $M$ is symmetric on $\D$.~\vspace{1mm}\\
{\bf ii)} Given $(x , y) \in \D$, since the corresponding
sequences ${(u_n)}_n$ and ${(v_n)}_n$ are respectively
non-decreasing and non-increasing and since $M(x , y)$ is the
common limit of ${(u_n)}_n$ and ${(v_n)}_n$ then we have:
$$u_0 \leq M(x , y) \leq v_0 ,$$
that is:
$$\min(x , y) \leq M(x , y) \leq \max(x , y)$$
{\bf iii)} Let $(x , y) \in \D$, fixed. Suppose that $M(x , y) =
x$ and show that $x = y$. Let us argue by contradiction; so assume
that $x \neq y$. Since $M_1$ and $M_2$ are mean functions on $\D$,
we have (according to the third axiom of mean functions):
\begin{equation}\label{eq7}
M_1(x , y) \neq x ~~\text{and}~~ M_2(x , y) \neq x
\end{equation}
We distinguish the following two cases:~\vspace{1mm}\\
\underline{1\textsuperscript{st} case:} (if $x < y$)\\
In this case, we have $M(x , y) = x = \min(x , y) = u_0$. So the
sequence ${(u_n)}_n$ is non-decreasing and converges to $u_0$. It
follows that ${(u_n)}_n$ is necessarily  constant and we have in
particular $u_1 = u_0$, that is:
$$\min(M_1(x , y) , M_2(x , y)) = x ,$$
which contradicts (\ref{eq7}).~\vspace{1mm}\\
\underline{2\textsuperscript{nd} case:} (if $x > y$)\\
In this case, we have $M(x , y) = x = \max(x , y) = v_0$. So the
sequence ${(v_n)}_n$ is non-increasing and converges to $v_0$. It
follows that ${(v_n)}_n$ is necessarily constant and we have in
particular $v_1 = v_0$, that is:
$$\max(M_1(x , y) , M_2(x , y)) = x ,$$
which again contradicts (\ref{eq7}).\\
Thus, we have $x = y$, as required.~\vspace{1mm}\\
From i), ii) and iii), we conclude that $M$ is effectively a mean function on $\D$.~\vspace{1mm}\\
Next, let us show that $M$ satisfies the functional equation
$M(M_1 , M_2) = M$. To do so, we constat that the fact to change
in the definition of the sequences ${(x_n)}_n$ and ${(y_n)}_n$ a
couple $(x , y)$ of $\D$ by the couple $(M_1(x , y) , M_2(x , y))$
just amounts to shift by one term these sequences (namely we
obtain ${(x_{n + 1})}_n$ instead of ${(x_n)}_n$ and ${(y_{n +
1})}_n$ instead of ${(y_n)}_n$). Consequently, this changing
conserves the common limit of the two sequences ${(x_n)}_n$ and
${(y_n)}_n$ (which is $M(x , y)$); that is:
$$M(M_1(x , y) , M_2(x , y)) = M(x , y) .$$
Since this last equation holds for all $(x , y) \in \D$, we have
$M(M_1 , M_2) = M$, as required.\\
It finally remains to show that $M$ is the unique mean on $\D$
satisfying the functional equation $M(M_1 , M_2) = M$. Let $M'$ be
a mean function on $\D$, satisfying $M'(M_1 , M_2) = M'$ and let
us show that $M'$ coincides with $M$. So, let $(x , y) \in \D$,
fixed and let us show that $M'(x , y) = M(x , y)$. We associate to
$(x , y)$ the sequence ${(x_n , y_n)}_{n \in \mathbb{N}}$ of
couples of $\D$, defined like in the second statement of the
theorem. Using the relation $M'(M_1 , M_2) = M'$, we have:
$$M'(x , y) = M'(x_1 , y_1) = M'(x_2 , y_2) = \cdots = M'(x_n , y_n) = \cdots$$
But since $M'$ is a mean on $\D$, it follows that for all $n \in
\mathbb{N}$, we have:
$$\min(x_n , y_n) \leq M'(x , y) \leq \max(x_n , y_n) .$$
Finally, since (according to the first part of this proof) $x_n$
and $y_n$ tend to $M(x , y)$ as $n$ tends to infinity, then by
tending $n$ to infinity in the three hand-sides of the previous
double inequality, we get:
$$M'(x , y) = M(x , y) ,$$
as required. This confirms the uniqueness of $M$ as a mean
function on $\D$ satisfying $M(M_1 , M_2) = M$. The proof is
complete.\penalty-20\null\hfill$\blacksquare$\par\medbreak

From Theorem \ref{t1}, we derive the following important
corollary:
\begin{coll}\label{coll2}
Let $\mathfrak{M}$ be a mean function on $\D = I^2$, with values
in $I$. Then, there exists a unique mean function on $\D$
satisfying the functional equation:
$$M\left(\frac{x + y}{2} , \mathfrak{M}(x , y)\right) = M(x , y) ~~~~~~ (\forall (x , y) \in \D) .$$
In addition, for all $(x , y) \in \D$, $M(x , y)$ is the common
limit of the two real sequences ${(x_n)}_n$ and ${(y_n)}_n$
defined by:
$$\left\{\!\begin{array}{l}
x_0 = x ~,~ y_0 = y \\
x_{n + 1} = \frac{x_n + y_n}{2}~~~~(\forall n \in \mathbb{N})\\
y_{n + 1} = \mathfrak{M}(x_n , y_n)~~~~(\forall n \in \mathbb{N})
\end{array}\right..$$
\end{coll}

\noindent{\bf Proof.} Since the metric space $(\M_\D , \d)$ is the
closed ball with center $\A$ and radius $1/2$ (see Proposition
\ref{prop5}), we have $\d(\mathfrak{M} , \A) \leq 1/2 < 1$. The
corollary then follows from Theorem \ref{t1} applied to the two
means $\A$ and $\mathfrak{M}$ on $\D$. The proof is
finished.\penalty-20\null\hfill$\blacksquare$\par\medbreak

\noindent{\bf Definition.} Given $\mathfrak{M}$ a mean on $\D$, we
suggest to call the mean $M$ given by Corollary \ref{coll2}: the
mean $\mathfrak{M}$-arithmetic. So the mean $\G$-arithmetic is
nothing else than the $\AGM$ mean.~\vspace{1mm}

\noindent{\bf Remark.} Since the metric space $(\M_\D , \d)$ is a
closed ball with radius $1/2$ then the distance between any two
means $M_1$ and $M_2$ on $\D$ is at most equal to $1$. It follows
that the unique case for which Theorem \ref{t1} cannot apply is
the extremal case where $\d(M_1 , M_2) = 1$ (remark also that in
such case the means $M_1$ and $M_2$ are obligatory on the border
of $\M_\D$). However, the functional middle between two means
$M_1$ and $M_2$ on $\D$ may exist and be unique even if the
condition $\d(M_1 , M_2) = 1$ is fulfilled. Indeed, taking $\D =
{(0 , + \infty)}^2$, $M_1 = \G$ and $M_2(x , y) = x + y - \sqrt{x
y}$ ($\forall (x , y) \in \D$), we easily verify that $\d(M_1 ,
M_2) = 1$ although the functional middle of $M_1$ and $M_2$ exists
and it is unique (it is simply the arithmetic mean $\A$ on
$\D$).~\vspace{1mm}

In the following theorem, we establish the existence and the
uniqueness of the functional middle of two mean functions $M_1$
and $M_2$ on $\D$ by changing the condition ``$\d(M_1 , M_2) <
1$'' (required by Theorem \ref{t1}) by the condition with
different nature which simply requires to $M_1$ and $M_2$ to be
continuous on $\D$.
\begin{thm}\label{t2}
Suppose that $I$ is an interval of $\R$ and let $M_1$ and $M_2$ be
two mean functions on $\D = I^2$ with values in $I$. We suppose
that $M_1$ and $M_2$ are continuous on $\D$. Then, there exists a
unique mean function $M$ on $\D$ satisfying the functional
equation:
$$M(M_1 , M_2) = M .$$
In addition, for all $(x , y) \in \D$, $M(x , y)$ is the common
limit of the two real sequences ${(x_n)}_n$ and ${(y_n)}_n$
defined by:
$$\left\{\!\begin{array}{l}
x_0 = x ~,~ y_0 = y \\
x_{n + 1} = M_1(x_n , y_n)~~~~(\forall n \in \mathbb{N})\\
y_{n + 1} = M_2(x_n , y_n)~~~~(\forall n \in \mathbb{N})
\end{array}\right..$$
\end{thm}

\noindent{\bf Proof.} Let $(x , y) \in \D$ fixed and let
${(x_n)}_n$ and ${(y_n)}_n$ be the two real sequences related to
$(x , y)$ which are introduced in the second part of the theorem.
Let also ${(u_n)}_n$ and ${(v_n)}_n$ be the two real sequences
defined by:
$$u_n := \min(x_n , y_n) ~~\text{and}~~ v_n := \max(x_n , y_n) ~~~~~~ (\forall n \in \mathbb{N}) .$$
It has been already shown during the proof of Theorem \ref{t1}
that ${(u_n)}_n$ is non-decreasing and that ${(v_n)}_n$ is
non-increasing. Since we have in addition $u_n \leq v_n$ $(\forall
n \in \mathbb{N})$ then ${(u_n)}_n$ bounded from above by $v_0$
and ${(v_n)}_n$ is bounded from below by $u_0$. It follows that
${(u_n)}_n$ and ${(v_n)}_n$ are convergent. Let $u = u(x , y)$ and
$v = v(x , y)$ denote the respective limits of ${(u_n)}_n$ and
${(v_n)}_n$ (so $u$ and $v$ lie in $[u_0 , v_0] =
[\min(x , y) , \max(x , y)] \subset I$). \\
Now, since $M_1$ and $M_2$ are symmetric on $\D$ (as mean
functions on $\D$), we have for all $n \in \mathbb{N}$:
$$x_{n + 1} = M_1(u_n , v_n) ~~\text{and}~~ y_{n + 1} = M_2(u_n , v_n) .$$
This implies (according to the hypothesis of continuity of $M_1$
and $M_2$ on $\D$) that the two sequences ${(x_n)}_n$ and
${(y_n)}_n$ are also convergent and that their respective limits
are $M_1(u , v)$ and $M_2(u , v)$.\\
Then, by tending $n$ to infinity in the two hand-sides of the
relation $x_{n + 1} = M_1(x_n , y_n)$, we obtain (according to the
continuity of $M_1$) that
$$M_1(u , v) = M_1\left(M_1(u , v) , M_2(u , v)\right) ,$$
which implies (according to the third axiom of mean functions)
that
$$M_1(u , v) = M_2(u , v) .$$
The two sequences ${(x_n)}_n$ and ${(y_n)}_n$ thus converge to a
same limit. Denoting (for all $(x , y) \in \D$) $M(x , y)$ the
common limit of ${(x_n)}_n$ and ${(y_n)}_n$, we show in the same
way as in the proof of Theorem \ref{t1} that $M$ is a mean
function on $\D$ and that it is the unique mean on $\D$ which
satisfies the functional equation $M(M_1 , M_2) = M$. The proof is
achieved.\penalty-20\null\hfill$\blacksquare$\par\medbreak

~\\


\begin{thebibliography}{9}
\bibliographystyle{plain}
\bibitem{bor}
J. M. Borwein, and P. B. Borwein, {\it Pi and the AGM, (A study in
Analytic Number Theory and Computational Complexity)}, John Wiley
\& Sons Inc., New York, 1987.
\end{thebibliography}
\end{document}